\documentclass[10pt,a4paper]{article}
\usepackage[utf8]{inputenc}
\usepackage{amsthm}
\usepackage{amsmath}
\usepackage{amsfonts}
\usepackage{amssymb}
\usepackage{url}
\usepackage{graphicx}
\usepackage{listings}
\usepackage{physics}
\usepackage{setspace}
\usepackage[marginratio=1:1]{geometry}

\newtheorem{remark}{Remark}[section]
\newtheorem{lemma}{Lemma}[section]
\newtheorem{corollary}{Corollary}[section]
\newtheorem{theorem}{Theorem}[section]

\newtheorem*{theorema}{Theorem A}
\newtheorem*{theorem*}{Theorem}
\newtheorem{definition}{Definition}

\author{Alexander Kushkuley \\ kushkuley@gmail.com}
\title{A Remark on Random Vectors and Irreducible Representations  }
\begin{document}
\maketitle

\begin{abstract}
\noindent The expectation of a squared scalar product of two   random independent unit  vectors that are uniformly distributed on a unit sphere in
$\mathbb{R}^n $ is  equal to $1/n$. We show  that this is a characteristic property of random unit vectors defined on  invariant probability subspaces  of irreducible representations of compact Lie groups.  We also discuss a relation of this  fact to some properties of random invariant tensors
\end{abstract}

\section{Introduction}
   \numberwithin{equation}{section}

 Let a compact Lie group $G$ linearly act on a real vector space $V$. We  will always  assume that  scalar product $ <,> $ on $V$  is $G$-invariant (cf. e.g \cite{Adams}) and we will identify elements of $G$ with linear operators on $ V $, denoting an action of $ g \in G $ on $v \in V $ by $g\circ v $ or just $gv$.

\bigskip\noindent Let $ \mathcal{U} \subset V $ be a $G$-invariant subset. Assume that there is a probability (Radon) measure $\mu$ defined on $\mathcal{U} $. The measure $\mu$ is called $G$-invariant if $ \mu (g\mathcal{U}') = \mu(\mathcal{U}') $ for any measurable  $\mathcal{U}' \subset \mathcal{U} $ and any $ g \in G $  (cf. e.g. \cite{Federer}, \cite{Leng}). Such a subset $ \mathcal{U}$ will be called an invariant probability subspace of $ V $.

 \bigskip\noindent
  We will sometimes omit mentioning of a group $G$ and just say that a measure $\mu$ and a set $ \mathcal{U} $ are invariant.
 A random vector $x$  distributed on $\mathcal{U} $ according to the invariant measure $\mu$ will thus  be called $G$-invariant random vector or just \emph{invariant random vector}.

\bigskip\noindent
 The expectation of a random vector $x$    defined on a probability space  $ \mathcal{U} \subset V $ will be denoted by $ \mathbb{E}_{\mu}(x) $ or  just $ \mathbb{E} (x) $. It is clear from definitions that for invariant random vectors  $ \mathbb{E}(gx) =  \mathbb{E}(x)  $.

\bigskip\noindent
 Let $ \chi $ be a normalized Haar measure on $G $  ( cf. e.g. \cite{Federer}).
An invariant projector (Reynolds operator) $ P^G : V \rightarrow V^G  $ onto the subspace of $G$-fixed points in $V$ can be  written as
\begin{equation}
P^G(v) = \int_G gv d\chi(g)  \nonumber
\end{equation}
for any $ v \in V $ (cf. e.g. \cite{Adams}, \cite{Collins2}).

\bigskip\noindent
The following statement is well known (cf. e.g. \cite{Collins1}-\cite{Collins2})

\begin{lemma}  If $x$ is a random vector defined on any invariant probability subspace of a real representation $V$ of a compact (Lie) group $G$, then

\begin{equation}
	\mathbb{E}(x) = \mathbb{E}( P^G (x) )  \nonumber
\end{equation}
\end{lemma}
\noindent Proof.  By invariance of $x$  and by linearity of expectation
\begin{equation}
\mathbb{E}( P^G (x) )  = \mathbb{E}  \left ( \int_G gx d\chi(g) \right)   =
\int_G \mathbb{E}(gx) d\chi(g) = \int_G \mathbb{E}(x) d\chi(g) =  \mathbb{E}(x)  \nonumber
\end{equation}
We mention some obvious corollaries
\begin{corollary}
	If a real representation of a compact Lie group has no non-trivial fixed points then the expectation of an invariant random vector defined on any invariant probability subspace of $V$ is zero.
\end{corollary}

\begin{corollary}
    If the representation of $G$ in $V$ is irreducible then  $
    \mathbb{E}( <x,y> ) =  0   \nonumber
    $ for  independent invariant  random vectors $x,y$ defined on any  invariant probability subspace  $ \mathcal{U} \subset S(V) $

\end{corollary}
\noindent To prove Corollary 1.2 note that
$    \mathbb{E}( <x,y> ) = \; <  \mathbb{E}(x),  \mathbb{E}(y) >
$
and apply Corollary 1.1

\subsection{Statement of the main result}
Let $G$ be a compact Lie group
and let $S(V)$ be a unit sphere in a vector space $V$ of a real $G$-representation of dimension $n$.

\begin{theorem}
\bigskip\noindent

\begin{itemize}

\item[(i)] For   independent invariant random vectors $x,y$ defined on any  invariant probability subspace  $ \mathcal{U} \subset S(V) $
		\begin{equation}
			\mathbb{E}( <x,y>^2 ) \geq 1/n   \nonumber
		\end{equation}
\item[(ii)] If the representation of $G$ in $V$ is irreducible then for  independent invariant  random vectors $x,y$ defined on any  invariant probability subspace  $ \mathcal{U} \subset S(V) $
	\begin{equation}
	\mathbb{E}( <x,y>^2 ) =  1/n   \nonumber
	\end{equation}

\item[(iii)] If there is a non-trivial linear invariant subspace  $W \subset V$ then  there is an invariant probability subspace  $ \mathcal{U} \subset S(V) $ and   independent invariant random vectors $x,y$ defined on $ \mathcal{U} $ such that
\begin{equation}
\mathbb{E}( <x,y>^2 ) >  1/n    \nonumber
\end{equation}
	\end{itemize}
 Conclusion: a real representation $V$ of a compact Lie group
 $G$ is irreducible if and only if
 \begin{equation}
 \mathbb{E}( <x,y>^2 ) =  1/n    \nonumber
 \end{equation}
  for any independent invariant random vectors $x,y$ defined on  any invariant probability subspace of the unit sphere

\end{theorem}
\noindent A proof of Theorem 1.1 and some  examples are provided in the next section. As we will see, the subject matter is naturally related to some properties of  tensor invariants. The last section (which is borrowed from \cite{On}) contains a brief discussion of this relation.
\begin{remark}
    The statement (ii) is well known for random vectors invariant with respect to the full orthogonal group, i.e. in the case of  $ G = O(n) $ and  $ \mathcal{U} = S(V)$ (cf. e.g. \cite{Al2}, \cite{Meckes}, Lemma 3.3 and Remark 3.1 below)
\end{remark}

\section{Proof of theorem 1.1}
We start with some well known   general observations (cf. e.g. \cite{Adams}).
\subsection{Preliminaries}

Let $V^{*} $ be a space of real linear functions on $ V $. This  is a space of a conjugate  $G$-representation, i.e. by definition
\begin{equation}
 gv^{*}(w) \; = \; v^{*}(g^{-1} w), \;v^{*} \in V^{*}, \; w \in V    \nonumber
\end{equation}

\bigskip\noindent
Consequently $G$ acts on $ V \otimes V^{*} $ and due to isomorphisms
\begin{equation}
V \otimes V^{*} \approx Hom(V,V) \approx Mat_n(V) \equiv M_n
\end{equation}
this action can be identified with the action of $G$ on algebra $ M_n = M_n( \mathbb{R})$ of  $n\times n $ real matrices, defined by the rule
  \begin{equation}
 g \circ a  = gag^{-1}, \; a \in M_n
  \end{equation}
Using  $G$-invariant scalar product, we can identify real vector spaces $V$ and $V^{*} $ (in orthogonal coordinates,   $v \in V$ is represented by a column vector and corresponding  $v^{*} \in V^{*} $ is represented by a transposed row vector $ v^{T}$).
We further notice that
for $ u_1, u_2 \in V $ and $  v_1, v_2 \in V^{*} $ we have
$ (u_1 \otimes v_1 )  (u_2 \otimes v_2 )  = v_1(u_2) u_1 \otimes v_2 = \; <v_1,u_2> u_1 \otimes v_2  $ and  it easy to see that both isomorphisms in (2.1) are isomorphisms of algebras.

\bigskip\noindent
The scalar product in $  V \otimes  V^{*}  \approx M_n$ is given by
\begin{equation}
< u_1 \otimes v_1 , u_2 \otimes v_2 > \; = \; <u_1,u_2><v_1, v_2 >
\end{equation}
and for matrices $ a, b \in M_n $ the formula (2.3) is equivalent to
\begin{equation}
<a,b> \; = \; Tr( ab^T )
\end{equation}
where $ Tr(.) $ stands for the trace of a matrix.
\noindent This scalar product is clearly $G$-invariant.
\bigskip\noindent

\noindent The following  lemma can be verified by a direct computation.
\begin{lemma}
	An orthogonal projection of a matrix $b \in M_n $ onto the one-dimensional  subspace spanned by a matrix $a \in M_n$ is given by
	\begin{equation}
		\frac{Tr(ab^T)}{Tr(aa^T)} \;a  \nonumber
	\end{equation}
	In particular, 	an orthogonal projection of a matrix $b$ onto the one-dimensional subspace spanned by a unity matrix $I_n$ is equal to
 	\begin{equation}
 	\frac{Tr(b)}{n} \; I_n     \nonumber
 	\end{equation}
\end{lemma}

\noindent There is a $G$-equivariant diagonal  map
\begin{equation}
\pi : V \rightarrow M_n, \;\; \pi( v ) = v \otimes v^{*}, \; v \in V   \nonumber
\end{equation}
that by (2.3)-(2.4) takes the unit sphere $ S(V) $ into a unit sphere in $    S(M_n) $.  If  $ \mathcal{U} \subset V $ is an invariant probability subspace in $ S(V) $ with probability measure $ \mu $ then its  image $ \pi (  \mathcal{U} ) $     is a $G$-invariant probabilty subspace in $  S(M_n)  $ with a $G$-invariant push-forward probability measure $ \hat{\mu} = \mu \circ \pi^{-1} $  (cf. \cite{Federer}). Thus we can say that
\begin{lemma}
	 If a random unit vector $x $ is invariant
	 then the random unit vector $ x \otimes x^{*} $  is also  invariant
\end{lemma}
\subsection{Expectation of $x\otimes x^{*} $ and conclusion of the proof }
Let $\mu$ be an invariant probability measure defined on $ \mathcal{U} \subset S(V) $. By (2.3), Lemma 2.1 and independence of $ x $ and $y $ we have
\begin{equation}
	\mathbb{E}_{\mu}( <x,y>^2 ) = 	\mathbb{E}_{  \hat{\mu}} (
	 <x \otimes x^{*}, \;  y \otimes y^{*} > )  \; = \;\;
	  < \mathbb{E}_{  \hat{\mu}} (x \otimes x^{*}), \; \mathbb{E}_{  \hat{\mu}} ( y \otimes  y^{*} ) >
\end{equation}
 where as above,  $ \hat{\mu} = \mu \circ \pi^{-1} $ is a push-forward probability measure. It follows from  (2.1) and (2.2)    that the unity matrix $I_n $ belongs to the linear subspace of fixed points $ M^{G}_{n} $. Let $ \mathcal{I} \subset  M^{G}_{n} $ be the one-dimensional   subspace of scalar matrices. We have an  orthogonal decomposition $  M^{G}_{n} = \mathcal{I} \oplus M' $ where  $ M' $ is some (possibly zero and obviously invariant) linear
subspace of  $ M^{G}_{n} $. It is clear  that an orthogonal (and hence invariant)  projection
of the random vector $x \otimes x^{*} $ onto  $ M^{G}_{n} $ is a sum of
its orthogonal projection  onto $ \mathcal{I} $ and its orthogonal    projection $ x' $ onto $ M' $.  Since $x$ is a unit vector, $ Tr( x \otimes x^{*} ) = 1 $ and applying lemma 2.2 we get
\begin{equation}
\mathbb{E}_{  \hat{\mu} } ( y \otimes y^{*} ) \;  = \; \mathbb{E}_{  \hat{\mu} } ( x \otimes x^{*} ) \;  = \; \mathbb{E}_{  \hat{\mu}}
\left( \frac{ Tr( x \otimes x^{*} ) } {n} \right)  I_n \; + \; \mathbb{E}_{  \hat{\mu}} (x') \; = \; \frac{I_n}{n}  \;  + \; z
\end{equation}
where $ z \in M' $. Statement (i) now follows from (2.5) and (2.6):
\begin{eqnarray}
	\mathbb{E}_{\mu}( <x,y>^2 ) \; = \;  < \frac{I_n}{n},  \frac{I_n}{n} > \;  + \;  <z,z> \; = \; 1/n \;\; + \; <z,z> \;\; \geq \; 1/n
\end{eqnarray}
Before going forward, observe that statement (i) and its proof are valid even for trivial $G$-representations. Below, we will make use of the  group-invariant proof of the statement (i). For now, we have in any case

\begin{corollary}
If $x,y$ are   independent  random unit vectors defined on any  probability subspace of  $ S(\mathbb{R}^n) $ then $\mathbb{E}( <x,y>^2 ) \geq 1/n  $.

\end{corollary}

\bigskip\noindent
To prove statement (ii), recall that if $V$ is an irreducible representation of $G$ then the linear subspace of fixed points  $ M^{G}_{n} $ of the $G$-action (2.2),  is isomorphic either to the  field of real numbers $ \mathbb{R} $ or to the complex field $ \mathbb{C} $ or to the field of quaternions $ \mathbb{H} $ (cf. e.g \cite{Adams}).
In the first case, the fixed point space is exactly  $ \mathcal{I} $ and (ii) follows from (2.7) where one can set $ z = 0 $. Something similar happens in other two cases.

We will consider only the case of the field of quaternions as it will become obvious that the case of complex numbers is covered by the quaternionic case. Hence, suppose that $ M^{G}_{n} $ is isomorphic to a fileld of quaternions  $ \mathbb{H} $ and therefore, there is a linear  basis of  $ M^{G}_{n} $  that consists of four real matrices $ I_n,  \mathfrak{i} ,  \mathfrak{j},  \mathfrak{k } \in M^{G}_{n} \subset  M_n $ that satisfy  relations of the algebra of quaternions, in particular  $  \mathfrak{i}^2  =  \mathfrak{j}^2  = \mathfrak{k }^2 = -I_n $.  It is not hard to figure out that matrices $  \mathfrak{i},  \mathfrak{j}, \mathfrak{k } $ are skew-symmetric and hence (cf. (2.4)) are orthogonal to $ I_n $ (and to each other, since e.g.
 $ Tr( \mathfrak{i} \mathfrak{j}^T) = -Tr( \mathfrak{k}) = 0 $).

On the other hand, the matrix  $x \otimes x^{*}$ is symmetric and therefore its projection onto a linear space of  skew-symmetric matrices is zero.  The conclusion is, that for irreducible representations, the extra term $z$ in (2.6)-(2.7) vanishes in all cases.

\bigskip\noindent
The proof of (iii) is very simple. Let $m, \; (0 < m < n) $ be a dimension of invariant subspace $W$ and let $ W'$ be an orthogonal complement of $W$ in $V$. Take, for example, an obviously invariant random vector $ x \in S(V) $ that has zero projection on $W'$ and is uniformly distributed on $ S(W) $. Then, by (i) the expectation of $ x $ is no less than $1/m > 1/n $.

\begin{corollary}
    Under  conditions of Theorem 1.1 (ii) \begin{equation}
        \mathbb{E} ( x \otimes x^{*} ) \;  = \; \mathbb{E}
        \left( \frac{ Tr ( x \otimes x^{*} ) } {n} \right)  I_n \; = \; \frac{I_n}{n}  \nonumber
    \end{equation}
\end{corollary}
\begin{corollary}
	Under conditions of Theorem 1.1 (ii), let Euclidean coordinates  of a random unit vector $x  \in \mathcal{U} \subset S(V) $ be $ x_1, x_2, \cdots,  x_n $. Then $ \mathbb{E}( x_i x_j ) = 0 $ if $ i \neq j $ and   $\mathbb{E}( x_i^{2}  ) = 1/n $ for all $ i, j  = 1,2, \cdots ,  n $.
\end{corollary}
\noindent The next corollary does not require  random vectors to be of unit length.
\begin{corollary}
    Let $ G $ be a compact group acting irreducibly on  $ \mathbb{R}^n$ and let $ \mathcal{U}, \mathcal{V} \subset \mathbb{R}^n $ be any $G$-invariant probability subspaces.   For any independent invariant random vectors $ x \in \mathcal{U}, \; y \in \mathcal{V} $
    	\begin{equation}
        \mathbb{E}( <x,y>^2 ) =  \frac{   \mathbb{E}(\norm{x}^2) \mathbb{E}(\norm{y}^2) }{n} \nonumber
    \end{equation}
\end{corollary}
\noindent The proof easily follows  from (2.5) and  Lemma 2.1 if one recalls that $ Tr( x \otimes x^{*} ) = \norm{x}^2 $
\subsection{Some Examples}

\subsubsection{Orthogonality relations for real irreducibe representations}
 Let $ \phi : G \rightarrow \textnormal{GL}(V_1) $
 and $ \psi : G \rightarrow \textnormal{GL}(V_2)  $ be two \textbf{complex} irreducible  representations (irreps) of our compact group $G $. Without loss of generality we will assume that both irreps $ \phi $ and $ \psi$ are unitary. Denote by
 $ \phi_{ij}(g), \;  i,j = 1, \cdots , \dim V_1 $ corresponding  matrix coefficients.  As complex functions on $G$  these coefficients  can be viewed as  random variables with respect to Haar measure on $ G $.
 The classical orthogonality relations theorem for matrix coefficients can be stated as follows
 \begin{theorema}("Grand orthogonality theorem" - see e.g. \cite{Serre}, \cite{Simon}, \cite{groupprops})
     \begin{enumerate}
         \item[(a)] if $ \phi  $ is not equivalent to $\psi$ then
         $ \mathbb{E}(\phi_{ij} \overline{\psi_{kl}}) = 0 $ for all index pairs $   i,j = 1, \cdots , \dim V_1 $ and $ k,\;l = 1, \cdots ,  \dim V_2$
         \item[(b)] $ \mathbb{E}(\phi_{ij} \overline{\phi_{kl}}) = 0 $ if $ i \neq k $ or $ j \neq l , \; 1 \leq i,j,k,l \leq \dim V_1 $
         \item[(c)] $ \mathbb{E}(\phi_{ij}  \overline{\phi_{ij}}) = 1/\dim V_1 ,    \; 1 \leq i,j \leq \dim V_1 $
     \end{enumerate}
\end{theorema}
\begin{remark}
    The orthogonality theorem is, of course, valid for non-unitary complex irreps. However, in a general case, the expectation   $ \mathbb{E}(\phi_{ij}(g) \overline{\psi_{kl}(g)}) $ should be replaced with
   $ \mathbb{E}\left(\phi_{ij}(g)\psi_{lk}(g^{-1})\right) $ (cf. e.g. \cite{Serre}). Same remark can be made about real orthogonal representations (see below)
\end{remark}
\noindent For real (orthogonal) representations, the statement (b) of Theorem A  does not hold. A simple (counter) example is an irreducible   representation $R$  of the circle group $ SO(2) $ by plane rotations. The Haar measure on $ SO(2)$ is  $ dt/(2\pi) $ where $ t $ is  an angular parameter on a unit circle, and
 \begin{align}
    \frac{1}{2\pi}\int_{SO(2)} R_{11}(t) R_{22}(t)dt = -   \frac{1}{2\pi}\int_{SO(2)} R_{12}(t) R_{21}(t)dt = 1/2
 \end{align}
 Nevertheless, the statement (c) of Theorem A is still valid for real (orthogonal) irreps. Moreover, as a corollary of Corollary 2.2 we have
 \begin{corollary}(Orthogonality theorem for real irreps).
 Let  $\phi, \psi$ be \textbf{real} orthogonal irreps.
 In notation of Theorem A
 \begin{enumerate}
     \item[(a)] if $ \phi  $ is not equivalent to $\psi$ then
     $ \mathbb{E}(\phi_{ij}\psi_{kl} ) = 0 $ for all index pairs $   i,j = 1, \cdots , \dim V_1 $ and $ k,l = 1 \cdots ,  \dim V_2$
     \item[(b)] if $j \neq l $ then  $\mathbb{E}(\phi_{ij} \phi_{il}) =   \mathbb{E}(\phi_{ji} \phi_{li}) = 0, \; 1 \leq i,j,l \leq \dim V_1 $
     \item[(c)] $ \mathbb{E}(\phi_{ij}^2) = 1/n $ where $ n = \dim V_1 $ and $ 1 \leq i,j \leq \dim V_1 $
 \end{enumerate}
 \end{corollary}
 \noindent   We will need below the following general remark.
  \begin{remark}
     Chose an orthonormal basis $ e_1, \cdots , e_n $ in $ V_1 $ so that  $ \phi_{ij}(g) = < ge_i, e_j > $.
   The real vector space $V_1$ is equipped with
     $G$-invariant scalar product (cf. the introduction). Using this scalar product we will identify spaces $ V_1$ and $ V_1^{\ast} $. Under this identification the invariant (matrix) $ I_n $ corresponds to a rank-two tensor $ \sum _{i=1}^n e_i \otimes e_i $ which will be denoted by the same symbol ($I_n$)
 \end{remark}

\noindent  Now, the statements (b) and (c) of Corollary 2.5 can be deduced from Corollary 2.2 as follows:
\begin{align}
  \mathbb{E}(\phi_{ij} \phi_{il}) =
  \int_G  < ge_i, \; e_j > <ge_i, \; e_l > dg \; =
  \int_G  < ge_i \otimes ge_i, \; e_j \otimes e_l
  >dg \; =  \nonumber \\  = \; < \mathbb{E} (ge_i \otimes ge_i, \; g \in G ), \; e_j \otimes e_l >  \; = \; \frac{1}{n}<I_n ,\;  e_j \otimes e_l  > \; = \; \frac{1}{n}\delta_{jl} \nonumber
\end{align}
\noindent and similarly
\begin{align}
    \mathbb{E}(\phi_{ji} \phi_{li}) =
    \int_G  < ge_j, \; e_i > <ge_l, \;  e_i > dg \; =  \int_G  < e_j, \; g^Te_i > <e_l, \; g^Te_i > dg \; =  \nonumber \\
    = \int_G  < e_j \otimes e_l, \; g^Te_i \otimes g^T e_i
    >dg \;  = \; \frac{1}{n}< e_j \otimes e_l, \; I_n  > \; = \; \frac{1}{n}\delta_{jl} \nonumber
\end{align}

\begin{remark}
    The statement (a) of Theorem A depends only on the first part of Schur Lemma (triviality of an intertwining operator for non-equivalent irreps, cf. e.g. \cite{Adams},\cite{Serre},\cite{Simon}). Therefore, the statement (a) obviously  holds for real  irreps as well.   The corollary 2.5 is applicable to all three kinds of real irreps (absolutely irreducible, complex and quaternionic, cf. e.g. ibid.).
     There could be other relations between matrix coefficients of real irreps that depend on irrep kind.  For example,  (2.8) does not hold for absolutely irreducible irreps.
\end{remark}

\subsubsection{An orbit as a probability space }
 We will see shortly that  for $G$-orbits, the statement (ii) of Theorem 1.1  is   a special case of the above mentioned orthogonality relations for matrix coefficients (cf. Corollary 2.5).
 Let $V$ be a real $n$-dimensional representation of a compact Lie group $ G $.  Let  $ v \in S(V) $ and let $ H = G_v \subset G $ be a stationary subgroup of $v$.  The normalized Haar (probability) measure   on $G$  induces  unique up to a scalar multiplier $G$-invariant measure on a (left) factor-space
$ G/H $ and hence on $G$-orbit   $ \mathcal{O}_v = Gv \subset S(V) $ (cf. e.g. \cite{Leng}).  Normalizing thus defined measure, one gets a probability measure $ \mu = \mu_{\mathcal{O}} $ on $ \mathcal{O} $ that will be called below  \emph{orbital (probability) measure}. For orbits of irreducible representations the statement (ii) of Theorem 1.1 can be read
as

\begin{corollary}
Let $ V $ be a real $n$-dimensional irreducible representation of a compact Lie group $ G $. Take an orbit of  $ \mathcal{O} = \mathcal{O}_v = Gv \subset S(V) $ of a unit vector $ v \in S(V) $ and let $\mu$ be the orbital probability measure on  $ \mathcal{O}  $. Then
\begin{equation}
\int_{ \mathcal{O}} \int_{ \mathcal{O}} <x, y>^2 d\mu(x) d\mu(y) = 1/n
\end{equation}
\end{corollary}
\noindent which leads to
\begin{corollary}
Under conditions of Corollary 2.6
\begin{equation}
\int_{ \mathcal{O} } <x, v>^2 d\mu(x) = 1/n
\end{equation}
\end{corollary}\bigskip\noindent
The proof is straightforward.
  Setting    $ y = g_y v,   $  for some (depending on $y$)  $ g_y \in G  $ in (2.9) and using invariance of the orbital measure, we get
\begin{eqnarray}
\frac{1}{n} \; = \; \int_{\mathcal{O}} \int_{\mathcal{O}} <x, y>^2 d\mu(x) d\mu(y ) \; =  \; \int_{\mathcal{O}} \int_{\mathcal{O}} < g^{-1}_y x, v>^2 d\mu(x) d\mu( y) \; =  \nonumber \\ = \int_{\mathcal{O}}d\mu(y ) \int_{\mathcal{O}} <x, v>^2 d\mu(x)  \; = \;
\int_{\mathcal{O}} <x, v>^2 d\mu(x) \nonumber
\end{eqnarray}
\bigskip\noindent

\begin{remark}
Any unit vector $v$ is an element of some orthogonal basis. Hence, (2.10) is again an orthogonality relation for  corresponding matrix coefficient $   <gv,v>, \; g \in G  $  (cf. Corollary 2.5 (c))
\end{remark}

\noindent
 It is probably worth mentioning a discrete version of Corollary 2.7.

\begin{corollary}
Let $ V $ be a space of a real irreducible exact $n$-dimensional  representation of a finite group $G$. Then for any $ v \in S(V) $
\begin{equation}
\sum_{k \in G} <kv, v >^2 \; = \;  \frac{|G|}{n}
\end{equation}
\end{corollary}
\noindent
Proof.  The Haar measure on a finite group $G$ is a point-mass measure that assigns equal probabilities $ 1/|G|, $ to all $ g \in G $, where $|G|$ denotes the order of the group. Hence, if a stationary subgroup $G_v$ of  a unit vector  $ v \in S(V) $ is trivial,  the corresponding  orbital measure $ \mu $ is given by   $   \mu( gv) = 1/|G| $ for all $ g \in G $. Since our representation is exact, the space $V$  contains  a \emph{principal orbit} of length $ |G| $, actually, the set of principal orbits is everywhere dense in $V$ (cf. e.g. \cite{Bredon}). If  a principal orbit of $G$ passes through   $v \in S(V) $ then (2.11) is precisely a discrete analogue  of Corollary 2.7. Moreover, any orbit is a limit of a sequence of principal orbits and therefore the same formula (2.11) is valid in general.
\begin{remark}
    It well may be that Theorem 1.1 itself is  a consequence of  orthogonality relations (cf. Corollary 2.5).  One can deduce (2.9) from (2.10) and then use an appropriate covering of an invariant probability space by tubular neighborhoods of principal orbits (at least for nice probability subspaces $ \mathcal{U} \subset S(V)$, cf. e.g. \cite{Bredon})
\end{remark}
\noindent We conclude this section with two concrete numerical examples.
\subsubsection{Cyclic soubgroups of SO(2) (cf. (2.8))}
For the cyclic group $C_n$ that rotates $ \mathbb{R}^2 $ by angles  $ 2\pi k / n, \; k = 0, 1, \cdots, n-1 $, the formula (2.11) reads as
\begin{equation}
    \sum_{k=0}^{n-1} \cos^2( \phi + 2\pi k /n ) = n/2  \nonumber
\end{equation}
where  $ \phi \in [ 0, 2\pi ] $ is arbitrary.
This is not surprising since
\begin{equation}
    \sum_{k=0}^{n-1} \cos^2( \phi + 2\pi k /n ) = \frac{n}{2} +  \frac{1}{2} \sum_{k=0}^{n-1} \cos( 2\phi + 4\pi k /n ) )   \nonumber
\end{equation}
and if $ n > 2 $ the last sum on the right is equal to zero, e.g. by Corollary 1. The condition $ n > 2 $ is not restrictive, since real two-dimensional representation of $ \mathbb{Z}_2$ is not irreducible. For one-dimensional irrep of $ \mathbb{Z}_2$ the formula (2.11) computes the variance ($=1$) of unbiased Bernoulli variable with values in $ \{1,-1\}$.
\subsubsection{Minimal irreducible representation of $ S_n $ }
 The symmetric group $ S_n $ acts on  $ \mathbb{R}^n $ by permuting Euclidean coordinates.
It is well known that the restriction of this $S_n$-action onto $(n-1)$-dimensional  linear subspace orthogonal to the vector $(1, \cdots, 1)$ is irreducible when $ n \geq 2$.
  Let
$ \cos( x,y)$ denote the cosine of an angle between vectors $ x, y \in \mathbb{R}^n $ and let,
as usual, $\parallel x \parallel$ be the Euclidean  norm of a vector $ x \in \mathbb{R}^n $. Interpretation of Corollary 2.8 in this setting yields
\begin{corollary}
	Let $ n \geq 2 $. If Euclidean coordinates  of a vector $ x \in  \mathbb{R}^n $ add up to zero, then
\begin{equation}
 \sum_{\sigma \in S_n} \cos^2( x, \sigma x)  \; \equiv \; \frac {1}{ \parallel x \parallel^4 } \sum_{\sigma \in S_n} < x, \sigma x >^2 \;  = \;  n!/( n-1)   \nonumber
\end{equation}
\end{corollary}

\section{Random Tensors and Veronese Surface}

The linear space $V^{\otimes,m}$   of rank-$m$ tensors  has a natural $G$-invariant dot-product (cf. (2.3)) defined by the rule :
\begin{equation}
    < x_1 \otimes \cdots
    \otimes x_m , \; y_1 \otimes \cdots
    \otimes y_m > \; = \; <x_1, y_1> \cdots   \nonumber
    <x_m, y_m>
\end{equation}
where $
x_1 \otimes \cdots
\otimes x_m , \; y_1 \otimes \cdots
\otimes y_m $ are any two decomposable tensors in   $V^{\otimes,m}$.
\newline\newline\noindent As in Lemma 2.2  we have a $G$-equivariant    map
\begin{equation}
    S(V) \ni	 x \xrightarrow{\nu} x ^{\otimes,m} = x \otimes \cdots \otimes x \in
    S(V^{\otimes,m}), \; m \geq 2   \nonumber
\end{equation}
 The image $\mathfrak{V}_m $ of  this map will be called rank-$m$ \emph{Veronese} surface (cf. e.g. \cite{Veronese}) and as in Lemma 2.2, an invariant
 probability subspace  $ \mathcal{U} \subset S(V)$ is mapped  unto invariant probability subspace of the Veronese surface    $ \nu(\mathcal{U})  \subset \mathfrak{V}_m $. We thus have (cf. e.g. Lemma 2.2) an
invariant random Veronese tensor  $ x ^{\otimes,m} \in \mathcal{U}^{\otimes,m}  \subset \mathfrak{V}_m $.
By now, the following statement should be obvious (cf. (2.5))
\begin{lemma}
    For   independent invariant random vectors $x,y$ defined on any  invariant probability subspace  $ \mathcal{U} \subset S(V) $
    \begin{equation}
        \mathbb{E}( <x,y>^m )  \; = \;  < \mathbb{E} ( x ^{\otimes,m} ) ,\;  \mathbb{E} ( y ^{\otimes,m} ) > \; = \;
        \mathbb{E} ( \parallel P^G(
        x ^{\otimes,m} ) \parallel^2 )
    \end{equation}
\end{lemma}

\noindent A natural question, therefore,  is how to compute the expectation of  $G$-invariant random Veronese tensor.
 When $m=2$, a simple answer to this question is that projection of $\mathfrak{V}_2 $  onto the space  of $G$-invariants is just one point, namely
 $ (1/n)\sum _{i=1}^n e_i \otimes e_i $ (cf. Corollary 2.2 and Remark 2.2).
 \noindent It is hard to believe that when  $ m > 2 $ a similar simple answer is possible for  arbitrary irreducible actions.  However, when acting group $G$ is the  orthogonal group $ O(n) $, the projection of Veronese surface onto the space of $G$-invariants  is again a single point and Corollary 2.2 can be nicely generalized. We will present some definitions before making a precise statement
 \subsection{$O(n)$-invariants}
 Fix  an orthonormal basis $ e_1, e_2, \cdots e_n $ in $V$ once and for all (cf. Remark 2.2)
 \begin{definition}
     Let $\mathcal{P} = \{J_1, J_2, \cdots, J_k \} $ be a partition of the set of indexes $ J_m = \{1,2, \cdots, m \} $.  For $ v_1, \cdots v_k \in V $ denote by
     \begin{equation}
         v_1(J_1) \otimes \cdots \otimes v_k(J_k)
     \end{equation}
     a decomposable tensor in $V^{\otimes,m}$ that has vector $v_1$  at indexes enumerated by $J_1$, vector
     $ v_2$ at positions enumerated by   $J_2$ and so on.
 \end{definition}
 \noindent For example, if $m = 5, \;
 J_1 = \{1,3,5\} $ and $ J_2 = \{2,4\} $ then $ u(J_1) \otimes v(J_2) $ denotes the tensor  $ u \otimes v \otimes u \otimes v\otimes u $.

 If $ m = 2k $ and all part sizes $|J_i|, \; i = 1, \; \cdots, k $ are equal, the partition $\mathcal{P}$ is called  \emph{pairing} (cf. e.g.  \cite{Collins1}-\cite{Collins2}, \cite{Good}).
 In case of $ m = 2$ there is just  one pairing $ J_0 = (1,2) $ and there is just one corresponding $G$-invariant tensor
 \begin{equation}
     I( J_0) = 	\sum_{i=1}^n e_i(J_0) =
     \sum_{i=1}^n e_i \otimes e_i \; \in V\otimes V  \nonumber
 \end{equation}
 In general, any pairing $\mathcal{P} = \{J_1, J_2, \cdots, J_k \} $ on $J_{m=2k}$ defines a corresponding \textit{standard} $G$-invariant tensor
 \begin{equation}
     I(\mathcal{P}) = 	\sum_{i_1,i_2, \cdots , i_k = 1}^{n}  e_{i_1}(J_{1})
     \otimes e_{i_2}(J_{2}) \otimes \cdots e_{i_k}(J_{k})  \nonumber
 \end{equation}
 Let $ \Pi_m$ denote the the set of all pairings on the set of indexes $J$. When $G$ is the full orthogonal group $O(V) = O(n) $  the classical First Fundamental Theorem (FFT) of  Invariant Theory  (see e.g. \cite{Proc}, \cite{Good} )  states that
 \begin{theorem*}\textbf{(FFT for Orthogonal Group, cf. \cite{Proc}-\cite{Lehrer})}.
     There are no non-trivial $O(V) $ invariants in   $V^{\otimes,m}$ if $ m $ is odd.
     Standard invariants $I(\mathcal{P}), \; \mathcal{P} \in \Pi_m $ span the linear space $ \mathcal{I}_m $ of $ O(V) $ invariants in $V^{\otimes,m}$ if $m$ is even.
 \end{theorem*}

 \noindent Any pairing, e.g. $ \mathcal{P}_1 =  (1,2)(3,4) \cdots ( 2k-1, 2k) $, can be viewed as an element of the permutation group $ S_m $ of the index set $J_{m}$ and  $S_m$ acts on the set of pairings $  \Pi_m \subset S_m $ by conjugation. If $ \sigma \in S_m $ then $ \sigma^{-1}  \mathcal{P}_1 \sigma =
 (\sigma(1),\sigma(2)) \cdots ( \sigma(2k-1), \sigma(2k) ) $ and  it is easy to see  that the stationary subgroup of $  \mathcal{P}_1 $  is $ T_k \mathcal{S}_k $, where $ T_k \approx \mathbb{Z}_2^{k} $ is generated by transpositions $ (1,2), \cdots , (2k-1, 2k) $ and $ \mathcal{S}_k \approx S_k $ is a permutation group of the set of  $k$ pairs $ (1,2), \cdots,   (2k-1, 2k) $. Clearly, the group  $ S_k $ normalizes
$ T_k $   and  $ \mathcal{S}_k \cap T_k = \{1\} $. Hence, $\Pi_m $   is an $S_m$ orbit of $ \mathcal{P}_1 $ of length $L =(2k)!/(2^k k!) = (2k-1)!! \equiv 1 \cdot 3 \cdots 2k - 1$ (cf. e.g. \cite{Good}). The following lemma is now obvious.

 \begin{lemma}
Let  $ \mathcal{G} = ( < I(\mathcal{P}_i,  I(\mathcal{P}_j) > ), \; i, j  = 1, \cdots , L ) $ be a Gram matrix of the ordered set of standard $O(m=2k) $ invariants.
 \begin{enumerate}
     \item[(i)]   Given two pairings $ \mathcal{P}_1,  \mathcal{P}_2 \in \Pi_m $
     \begin{equation}
         < I(\sigma \mathcal{P}_1),  I(\mathcal{P}_2) > \; = \; < I(\mathcal{P}_1),  I( \sigma^{-1}\mathcal{P}_2) >  \nonumber
     \end{equation}
     for any $ \sigma \in S_m $
     \item[(ii)] Sums of elements of every row of $ \mathcal{G} $ are equal to each other
 \end{enumerate}

 \end{lemma}

\subsubsection{Projection of Veronese surface onto the space $O(n)$-invariants}
\noindent From now on, unless explicitly stated otherwise, we assume that the group $G$ is a full orthogonal group $O(V) \approx O(n)$.
For obvious reasons, everything that was said above remains valid in this special case.
Note that $S(V) $ is an orbit of $O(V)$ and  hence the $O(V)$-invariant probability measure on a unit sphere $S(V)$ is uniquely inherited from the Haar measure on $
O(n)$ (cf. e.g. \cite{Leng} and section 2.3.2).
By direct computation (see  for example \cite{Folland}, \cite{Meckes}) one gets
\begin{lemma} (\cite{Meckes}, \cite{Folland})
    For  any independent $O(n)$-invariant random vectors $x,y$ defined on   $ S(V) $ and any integer $ k \geq 1 $
    \begin{equation}
        \mathbb{E}\left( <x,y>^{2k} \right) = \frac{(2k-1)!!}{(n + 2k - 2)(n + 2k - 4) \cdots n}
    \end{equation}
\end{lemma}
\begin{remark}
    Note that for $ k = 1$ this is the same formula as in Theorem 1.1 (ii). The difference is, that in general, Lemma 3.3 requires invariance with respect to the full orthogonal group.
\end{remark}
\noindent Let $m = 2k $ as above. Denote the denominator of the right hand side of (3.3) by $ P(n,k ) $ and  recall that $ \mathcal{G} $ (cf. Lemma 3.1)
denotes the Gram matrix of the ordered standard set of invariants   $  I(\mathcal{P}_i), \;  i = 1, \cdots, L   $.
Let
\begin{equation}
    A_m =  \frac{1}{L} \sum_{i=1}^L  I(\mathcal{P}_i)  \in \mathcal{I}_m   \nonumber
\end{equation}
be the average of the set of all standard  invariants of $O(n) $.
\begin{theorem} ($m=2k$)

    \begin{enumerate}
        \item[(i)]  	The orthogonal projection of the Veronese surface $ \mathfrak{V}_m $ onto the linear space of invariants $ \mathcal{I}_m $ is a single point $ 	\frac{L}{P(n,k)}A_m$.
        \item[(ii)] Sum of the elements of any row of $\mathcal{G} $ is equal to $P(n,k)  $
    \end{enumerate}
      In other words, the  expectation of $O(n)$-invariant  random Veronese tensor $  x ^{\otimes,2k} $ is an average of standard $O(n)$-invariants  divided by the average of all entries of their  Gram matrix

\end{theorem}
\noindent Proof.  Since $ \mathfrak{V}_m $ is an orbit of $O(V) $, its orthogonal projection onto $\mathcal{I}_m$ is a  single point. To find this point, let
\begin{equation}
    P = \sum_{ i=1}^L \alpha_i  I(\mathcal{P}_i)  \in \mathcal{I}_m, \; \alpha_i \in \mathbb{R}   \nonumber
\end{equation}
be a projection of some Veronese tensor $ a ^{\otimes,m}, \; a \in S(V) $ onto
the space of invariants $\mathcal{I}_m$. The coefficients  $  \alpha_i, \; i = 1, \cdots , L  $ must satisfy the system of "normal equations"
\begin{equation}
    <P -  a ^{\otimes,m}, \;  I(\mathcal{P}_j) > \; = \; 0 , \; j = 1, \cdots , L
\end{equation}
\noindent It is obvious that  $<a ^{\otimes,m}, \;  I(\mathcal{P}_j )> \; =  1 $ for all $ j $ and therefore the system of  linear equations (3.4) can be rewritten as
\begin{equation}
    \mathcal{G}  \alpha = 1_{L}
\end{equation}
where $ \alpha$ is an unknown vector with real coordinates
$\alpha_i, \; i = 1, \cdots, L $ and
$ 1_L$ is the vector with all its $L$ coordinates equal to one. Denoting the unique row-sum of $ \mathcal{G}$ (Lemma 3.2 (ii)) by $N$ we see that the vector $(1/N) 1_L $
is a solution of the  system of equations (3.5).
The  value of $N$ can be  found by comparing (3.1) and (3.3).
Indeed, from the uniqueness of  point $P$ we get
\begin{align}
    P = 1/N \sum_{i=1}^L I(\mathcal{P}_i) = 	\mathbb{E} \left(  x ^{\otimes,m=2k} \right)
\end{align}

\noindent now by Lemma 3.1 (with $m=2k$) and  (3.6), (3.3)
\begin{align}
\text{"right hand side of (3.3)"} \equiv    L/P(n,k) \; = \;  \frac{1}{N^2}<\sum_{ \mathcal{P} \in \Pi_m} I(\mathcal{P}) , \sum_{ \mathcal{P} \in \Pi_m} I(\mathcal{P}) > \; = \; 	LN/N^2   \nonumber
\end{align}
Hence, we get $ N= P(n,k)$ finishing the proof.
\begin{remark} If $ m = 2k = 2 $ then $ L = 1 $, $P(n,k) = n $ and we see that for the full orthogonal group, Corollary 2.2 is a particular case of Theorem 3.1
\end{remark}
\begin{remark}
     The Gram matrix of the standard set of $O(n)$-invariants is degenerate for large $m$ (see e.g. \cite{Proc} or a nice summary in \cite{Lehrer}) and the linear system (3.5) could be overdetermined. This is not an obstacle, since we have a unique solution handy.
\end{remark}
\begin{remark}
    The Gram matrix $ \mathcal{G} $ plays an important role in Weingarten calculus (cf. e.g. \cite{Collins1}, \cite{Collins2}). Some properties of this matrix are analyzed in \cite{On})
\end{remark}

\end{document}